\documentclass[a4paper,12pt]{amsart}
\usepackage{graphicx,amssymb}
\usepackage{txfonts}
\usepackage{mathrsfs}
\usepackage{amssymb,amsmath,amsthm,color}
\usepackage{graphicx,mcite}
\usepackage{hyperref}
\usepackage{url}
\usepackage{setspace}
\newtheorem{theorem}{Theorem}[section]
\newtheorem{proposition}[theorem]{Proposition}
\newtheorem{lemma}[theorem]{Lemma}
\newtheorem{definition}[theorem]{Definition}

\newtheorem{remark}[theorem]{Remark}
\newtheorem{corollary}[theorem]{Corollary}

\def\be#1 {\begin{equation} \label{#1}}
\newcommand{\ee}{\end{equation}}

\def\sqw{\hbox{\rlap{\leavevmode\raise.3ex\hbox{$\sqcap$}}$%
\sqcup$}}
\def\findem{\ifmmode\sqw\else{\ifhmode\unskip\fi\nobreak\hfil
\penalty50\hskip1em\null\nobreak\hfil\sqw
\parfillskip=0pt\finalhyphendemerits=0\endgraf}\fi}

\newcommand{\R}{{\mathbb {R}}}

\newcommand{\N}{{\mathbb N}}
\newcommand{\Z}{{\mathbb Z}}

\newcommand{\T}{{\mathbb T}}

\author{K. Jotsaroop }
\email{jyotsy.12@gmail.com}
\address{Department of mathematics, Indian Institute Science Education and Research Mohali, India}

\author{Saurabh Shrivastava}
\email{saurabhk@iiserb.ac.in} 
\address{
Department of Mathematics\\
Indian Institute Science Education and Research Bhopal\\
Bhopal-462066, India}


\begin{document}

\title[Unimodular bilinear Fourier multipliers on $L^p$ spaces]{Unimodular bilinear Fourier multipliers on $L^p$ spaces}

\begin{abstract}
In this paper we investigate the boundedness properties of bilinear multiplier operators associated with unimodular functions of the form $m(\xi,\eta)=e^{i \phi(\xi-\eta)}$. We prove that if $\phi$ is a $C^1(\R^n)$ real-valued non-linear function, then for all exponents $p,q,r$ lying outside the local $L^2-$range and satisfying the H\"{o}lder's condition $\frac{1}{p}+\frac{1}{q}=\frac{1}{r}$, the bilinear multiplier norm
$$\|e^{i\lambda \phi(\xi-\eta)}\|_{\mathcal M_{p,q,r}(\R^n)}\rightarrow \infty,~ \lambda \in \R,~ |\lambda|\rightarrow \infty.$$
For exponents in the local $L^2-$range, we give examples of unimodular functions of the form $e^{i\phi(\xi-\eta)}$, which do not give rise to bilinear multipliers. Further, we also discuss the essential continuity property of bilinear multipliers for exponents outside local $L^2-$ range.
\end{abstract}

\maketitle
\section{Introduction \& preliminaries}
\subsection{Fourier multipliers} For $m\in L^{\infty}(\R^n)$ consider the linear operator  
\begin{eqnarray}\label{def10}
\widehat {S_m f} :=m \hat f,~f\in L^2(\R^n) \cap L^p(\R^n),
\end{eqnarray}
where $\hat{f}$ denotes the Fourier transform of $f$. 

The function $m\in L^{\infty}(\R^n)$ is called an $L^p-$Fourier multiplier if the associated linear operator $S_m$ is bounded on $L^p(\R^n)$. The space of all $L^p-$Fourier multipliers is denoted by $M_p(\R^n).$ It forms a Banach algebra with respect to the standard pointwise multiplication of functions and the norm is given by 
$$\|m\|_{M_p(\R^n)}:= \|S_m\|_{L^p \rightarrow L^p}.$$

The notion of Fourier multipliers on groups $\T^n$ and $\Z^n$ are defined similarly. 

It is an easy consequence of the Plancherel theorem that the space $M_2(\R^n)\cong L^{\infty}(\R^n)$. Further, it is also known that the space $M_1(\R^n)$ coincides with the algebra $A(\R^n)$ of functions $f$ of the form $f=\hat \mu$, where $\mu$ is a bounded Borel measure on $\R^n$. When $p\neq 1,2$, no such characterization of spaces $M_p(\R^n)$ is known. 

For a  measurable function $\phi:\R^n\rightarrow \R$ and $\lambda \in \R$ consider the unimodular function $e^{i \lambda \phi}.$ Such functions are of special interest in the theory of Fourier multipliers. The boundedness properties of the Fourier multiplier operators associated with functions $\phi(\xi)=|\xi|^{\alpha}$ play important role in PDEs. In particular, the cases $\alpha=1$ and $\alpha=2$ occur in the study of time evolution of the wave equation and the free Schr\"{o}dinger operator respectively. Another motivation to study the boundedness of such Fourier multiplier operators comes from the Beurling-Helson theorem, which raises the question of describing real-valued functions $\phi$ so that 
$$\|e^{i \lambda \phi}\|_{M_p(\R^n)}=\huge{O}(1),~~\lambda \in \R.$$
We refer to~\cite{h,o1,o2} and the references therein for an elaborate discussion about this point. 

L. H\"{o}rmander~\cite{h} proved that if $\phi$ is a $C^2(\R^n)$ real-valued function and 

\begin{eqnarray} \label{hor1} \|e^{i \lambda \phi}\|_{M_p(\R^n)}=\huge{O}(1),~~\lambda \in \R,
\end{eqnarray}
where $1<p<\infty, p\neq 2,$ then $\phi$ is a linear function. In fact, he showed that the same assertion holds true with the weaker assumption that  
\begin{eqnarray} \label{hor2} \|e^{i \lambda_m \phi}\|_{M_p(\R^n)}=\huge{O}(1),
\end{eqnarray}
where $\{\lambda_m\}$ is an unbounded sequence of real numbers. 

In \cite{h}, L. H\"{o}rmander conjectured that the above result holds for $\phi\in C^1(\R^n)$ as well. In 1994 V. Lebedev and A. Olevskii \cite{o1} settled this conjecture and proved that if $\phi\in C^1(\R^n)$ satisfies (\ref{hor1}) then $\phi$ is a linear function. This together with L. H\"{o}rmander \cite{h} proves that the condition $\phi\in C^1(\R^n)$ is sharp. Further, in~\cite{o2} V. Lebedev and A. Olevskii obtained the following generalization of this result: If 
$\phi:\R^n \rightarrow [0,2\pi)$ is a measurable function and 

\begin{eqnarray} \label{hor3} \|e^{i \lambda \phi}\|_{M_p(\R^n)}=\huge{O}(1),~~\lambda \in \mathbb Z,
\end{eqnarray}
where $1<p<\infty, p\neq 2,$ then $\phi$ is a linear function on domains complementary to some closed set, depending on $\phi$, of measure zero. Moreover, $\nabla \phi$ takes only finitely many values.

As a consequence of the above mentioned results along with the standard dilation argument we get that the function $e^{i  \phi}$ is not an $L^p$ Fourier multiplier, $1<p<\infty,~p\neq 2,$ if $\phi$ is  $C^1(\R^n)$ homogeneous non-linear function. 
At this point we would like to refer to~\cite{ben} for some positive results about boundedness of Fourier multipliers $e^{i|\xi|^{\alpha}}$ on modulation spaces. 

In this paper, we investigate the boundedness of bilinear multiplier operators associated with functions of the form $e^{i \phi}$ and provide answers to various questions concerning bilinear multipliers under consideration. 
\subsection{Bilinear multipliers}  \label{bil}
 Let $m(\xi,\eta)$ be a bounded measurable function on $\R^n\times \R^n$ and $(p,q,r),~0<p,q,r\leq \infty$ be a triplet of exponents. In what follows, we will always assume that $p,q,r$ satisfy the H\"{o}lder condition i.e. $\frac{1}{p}+\frac{1}{q}=\frac{1}{r}$. Consider the bilinear operator $M_m$ initially defined for functions $f$ and $g$ in a suitable dense class by  
\begin{eqnarray}\label{def1}
M_{m}(f,g)(x)= \int_{\R^n}\int_{\R^n}\hat{f}(\xi)\hat{g}(\eta) m(\xi,\eta)e^{2\pi i x\cdot(\xi+\eta)}d\xi d\eta.
\end{eqnarray}

We say that $M_m$ is a bilinear multiplier operator for the triplet $(p,q,r)$ if $M_m$ extends to a bounded operator 
from $L^p(\R^n)\times L^q(\R^n)$ into $L^r(\R^n),$ i.e. there exists a constant $C>0,$ independent of functions $f$ and $g,$ such that  
\begin{eqnarray*}\label{bm}
\|M_m(f,g)\|_{L^r(\R^n)}\leq C \|f\|_{L^p(\R^n)} \|g\|_{L^q(\R^n)}.
\end{eqnarray*}

The bounded function $m$ is said to be a bilinear multiplier symbol for the triplet $(p,q,r)$ if the corresponding operator $M_m$ is a bilinear multiplier operator for $(p,q,r).$ 

We denote by $\mathcal{M}_{p,q,r}(\R^n)$ the space of all bilinear multiplier symbols for the triplet $(p,q,r).$ Further, the norm of $m\in \mathcal{M}_{p,q,r}(\R^n)$ is defined to be the norm of the corresponding bilinear multiplier operator $M_m$ from $L^p(\R^n)\times L^q(\R^n)$ into $L^r(\R^n),$ i.e.
\begin{eqnarray*}
\|m\|_{\mathcal{M}_{p,q,r}(\R^n)}=\|M_m\|_{L^p(\R^n)\times L^q(\R^n)\rightarrow L^r(\R^n)}.
\end{eqnarray*}
For exponents $p,q,r\in [1,\infty]$, the real duality helps us consider two adjoint operators associated with a given bounded bilinear operator $M_m$ from $L^p(\R^n)\times L^q(\R^n)$ into $L^r(\R^n)$. These adjoint operators are given by
$$ \langle M_m^{*,1}(h,g),f\rangle: =\langle M_m(f,g),h\rangle ~\text{and}~ \langle M_m^{*,2}(f,h),g\rangle: =\langle M_m(f,g),h\rangle.$$

Note that $M_m^{*,1}$ is bounded from $L^{r'}(\R^n)\times L^q(\R^n)$ into $L^{p'}(\R^n),$ and  $M_m^{*,2}$ is bounded from $L^p(\R^n)\times L^{r'}(\R^n)$ into $L^{q'}(\R^n).$ Moreover, for bilinear multipliers $m(\xi,\eta)\in \mathcal{M}_{p,q,r}(\R^n)$, using the adjoint operators one can easily conclude that  $m(\xi+\eta,-\eta)\in \mathcal{M}_{r',q,p'}(\R^n)$ and $m(-\xi,\eta+\xi)\in \mathcal{M}_{p,r',q'}(\R^n)$ with norm exactly the same as that of $m$. 

The bilinear multipliers on the Torus group $\T^n$ and discrete group $\Z^n$ are defined similarly. The space of bilinear multipliers on $\T^n$ and $\Z^n$ will be denoted by $\mathcal{M}_{p,q,r}(\T^n)$ and $\mathcal{M}_{p,q,r}(\Z^n)$ respectively.

\subsection{Known results in bilinear setting}
Let $\phi$ be a measurable function defined on $\R^n$. For $f,g \in \mathcal S(\R^n)$ consider the bilinear operator
\begin{eqnarray}\label{bo}
B(f,g)(x):=\int_{\R^n}\int_{\R^n} \hat f(\xi)\hat g(\eta) e^{i\phi(\xi,\eta)}e^{2\pi i x\cdot(\xi+\eta)} d\xi d\eta
\end{eqnarray}
where $\phi$ is a non-linear function. 
These types of bilinear multipliers arise when we study the solution of the non-linear PDE
\begin{eqnarray*}
\partial_t u(t,x) + P(D)u(t,x) &=& |u(t,x)|^2\\
u(0,x)&=& f(x),
\end{eqnarray*}
where $P(D)$ is a quadratic function of $D=(\partial_{x_1},\partial_{x_2},\cdots,\partial_{x_n})$.
The solution of the above PDE is given by 
\begin{eqnarray*}
u(t,x)= e^{i tP(D)}f(x) + \int_{0}^{t}e^{i(s-t)P(D)}(u(s,.),\overline{u(s,.)} )(x)ds.
\end{eqnarray*}
It is therefore natural to study the boundedness properties of bilinear operators of the type defined above in (\ref{bo}).
 
We refer to the work of F. Bernicot and P. Germain in ~\cite{BG} for more details. In this paper, the authors established several results about the boundedness of bilinear oscillatory integral operators. In order to describe the relevant results to our context, we need the following definition. 
\begin{definition}(Local $L^2-$ range) The sets of exponents $(p,q,r)$ satisfying $2\leq p, q, r'\leq \infty$ and the H\"{o}lder condition $\frac{1}{p}+\frac{1}{q}=\frac{1}{r}$ is referred to as the local $L^2-$ range of exponents in the context of bilinear multipliers. We shall use the notation $L$ for this set. 
\end{definition}

For the bilinear multiplier operator $T_{\lambda}$ defined below, and the exponents $(p,q,r)$ in local $L^2-$ range,  
F. Bernicot and P. Germain proved the following theorem for bilinear oscillatory integral operators with symbols defined on $\R^2$ . 
\begin{theorem}\cite{BG}
Let us assume that
$$\partial_{\eta}\partial_{\xi}\phi\neq 0,$$
$$(\partial_{\eta}^2 -\partial_{\eta}\partial_{\xi})\phi\neq 0$$ and $$(\partial_{\xi}^2 -\partial_{\eta}\partial_{\xi})\phi\neq 0.$$
Then the bilinear oscillatory integral operator 
$$T_{\lambda}(f,g)(x):= \frac{1}{(2\pi)^{1/2}}\int_{\R}\int_{\R}\hat f(\xi)\hat g(\eta) e^{i\lambda\phi(\xi,\eta)}m(\eta,\xi)e^{2\pi i x\cdot(\xi+\eta)} d\xi d\eta$$
satisfies the following boundedness: for all exponents $(p,q,r)$ in the local $L^2-$ range, 
there exists a constant $C=C(p,q,r,\phi,m)$ such that for all $\lambda\neq 0$
$$\|T_{\lambda}(f,g)\|_{L^{r'}}\leq C |\lambda|^{-\frac{1}{2}}\|f\|_{L^{p'}}\|g\|_{L^{q'}}.$$

\end{theorem}

Further in this paper they study  variants of the  bilinear oscillatory integral operator of the form $T_{\lambda}(f,g)$ in different settings.
Note that the above is a bilinear analogue of dispersive estimates of the oscillatory integral operators in the linear case. However in this paper we are interested in studying the bilinear analogue of the linear $L^p$ multiplier of the form $e^{i\phi}$.
Motivated from the work of L. H\"{o}rmander~\cite{h}; V.~Lebedev and A. Olevskii~\cite{o1,o2} we study the boundedness of the operator 
\begin{eqnarray*}
T_{\phi,\lambda}(f,g)(x):=\int_{\R^n}\int_{\R^n} \hat f(\xi)\hat g(\eta) e^{i\lambda \phi(\xi-\eta)}e^{2\pi i x\cdot(\xi+\eta)} d\xi d\eta
\end{eqnarray*} from $L^p(\R^n)\times L^q(\R^n)$ into $L^r(\R^n)$.
In this paper we are mainly interested in bilinear multipliers of the type  $m(\xi,\eta)=\psi(\xi-\eta),$ where $\psi \in L^{\infty}(\R^n).$ 
The bilinear multipliers of this form satisfy a modulation invariant property and are of  particular interest in the theory of bilinear multipliers. The bilinear Hilbert transform is the key example of bilinear multiplier operators having modulation invariance. We refer to the celebrated works of M. Lacey and C. Thiele~\cite{lt1,lt2} for precise details on the bilinear Hilbert transform. 
In particular, we prove that if $\phi$ is a $C^1(\R^n)$ real-valued function satisfying the bilinear analogue of (\ref{hor1}) (see Theorem~\ref{mt}) for exponents lying outside local $L^2-$ range then $\phi$ is a linear function. The method of the proof of Theorem~\ref{mt} is based on the ideas by V.~Lebedev and A. Olevskii~\cite{o1,o2}. For exponents in local $L^2-$ range, we will give examples of non-linear functions $\phi$ for which the functions of the forms $e^{i \phi(\xi-\eta)}$ and $e^{i \phi(\xi,\eta)}$ do not give rise to bilinear multipliers. 

Further, in Section~\ref{ec} we establish essential continuity property of certain bilinears for exponents outside the local $L^2-$ range. 
\begin{remark} Unlike the linear case,  in the bilinear setting there does not exist a consistent range of exponents $p,q$ and $r$ for which the class of functions under consideration always gives rise to bilinear multipliers. 
\end{remark}
\begin{remark} Let $\phi_1$ and $\phi_2$ be linear and non-linear functions respectively. The H\"{o}lder's inequality yields that the function of the form $\phi(\xi,\eta)=\phi_1(\xi)\phi_2(\eta)$ gives rise to unimodular bilinear multiplier, in the sense as above, from $L^p(\R^n) \times L^2(\R^n)$ into $L^q(\R^n)$ for exponents $p$ and $q$ satisfying the H\"{o}lder condition $\frac{1}{p}+\frac{1}{2}=\frac{1}{q}$. The bilinear multipliers defined using functions of the form $\phi(\xi-\eta)$ cannot be factored in a similar fashion and need to be investigated separately.  
\end{remark}
\section{Main result and proofs}
The main result of this article is the following.
\begin{theorem}\label{mt} Let $(p,q,r)$ be a triplet of exponents outside the local $L^2-$range. If $\phi$ is a $C^1(\R^n)$ real-valued non-linear function, then
$$\|e^{i\lambda \phi(\xi-\eta)}\|_{\mathcal M_{p,q,r}(\R^n)}\rightarrow \infty,~ \lambda \in \R,~ |\lambda|\rightarrow \infty.$$
\end{theorem}

The theorem above can be deduced from the following lemma. 
\begin{lemma}\label{ml}
Let $\phi:\R^n\rightarrow \R$ be a measurable function. Suppose that there are $N$ cubes $Q_k\subset \R^n, k=1,2,\dots,N$ such that $\phi(t)=\langle \alpha_k, t\rangle +\beta_k$ for almost all $t\in Q_k$, the vectors $\alpha_k, k=1,2,\dots, N$ are all distinct and $\beta_k\in \R^n$. Then for any unbounded sequence of real numbers $\{\lambda_k\}_{k\in \mathbb N}$ we have 
\begin{eqnarray}\label{norm} \sup_{k\in \mathbb N} \|e^{i\lambda_k \phi(\xi-\eta)}\|_{\mathcal M_{p,q,r}(\R^n)} \geq N^{\gamma-\frac{1}{2}},
\end{eqnarray}
where $\gamma=\max\{\frac{1}{p},\frac{1}{q},\frac{1}{r'}\}$.
\end{lemma}

In order to prove the results stated above we would need the following auxiliary result about vector-valued bilinear operators.  
\begin{proposition} \label{vector} Let $\alpha=(\alpha_1,\alpha_2,\dots,\alpha_N)$ be an $N-$tuple of distinct vectors in $\R^n$ and $\rho>0$ be a positive real number. For vector-valued function $f=(f_1,f_2,\dots,f_N)$ and $g=(g_1,g_2,\dots,g_N)$ consider the bilinear operator 
\begin{eqnarray*}
S_{\alpha,\rho}(f,g)(x)=(f_1(\cdot+\rho\alpha_1)g_1(\cdot-\alpha_1), f_2(\cdot+\rho\alpha_2)g_2(\cdot-\alpha_2), \dots, f_N(\cdot+\rho\alpha_1)g_N(\cdot-\alpha_N))(x).
\end{eqnarray*}
Then the norm of the bilinear operator $S_{\alpha,\rho}$ satisfies the following 
$$\|S_{\alpha,\rho}\|_{p,q,r} \geq \max\{N^{1/p-1/2}, N^{1/q-1/2}\}.$$
\end{proposition}

{\bf Proof:}
Since $\alpha_k$'s are distinct and finitely many in number, we can choose a $\delta>0$ such that 
$$\rho |\alpha_j-\alpha_i|>\delta~~~\forall 1\leq i\neq j\leq N.$$
Choose $f=(f_1,f_2,\dots,f_N)$ and $g=(g_1,g_2,\dots,g_N)$ with $f_k=\chi_{_{B(0,\delta)}}$ and $g_k=\chi_{_{B((1+\rho)\alpha_k,\delta)}}$. Observe that $\|f\|_{L^p(l_2)}=N^{1/2} |B(0,\delta)|^{1/p}$ and  $\|g\|_{L^q(l_2)}=(N|B(0,\delta)|)^{1/q}$. 

Further, note that 
$f_k(\cdot+\rho \alpha_k)=\chi_{_{B(\rho\alpha_k,\delta)}}(\cdot)$ and $g_k(\cdot-\alpha_k)=\chi_{_{B(\rho \alpha_k,\delta)}}(\cdot)$.
It is easy to verify that 
\begin{eqnarray*}
\left( \sum\limits_{k=1}^N |f_k(x+\rho \alpha_k)g_k(x-\alpha_k)|^2\right)^{\frac{1}{2}}=\chi_{_{\bigcup_{k=1}^N B(\rho \alpha_k,\delta)}}(x)
\end{eqnarray*}
Therefore, 
$$\|S_{\alpha,\rho}(f,g)\|_{L^r(l_2)}=(N|B(0,\delta)|)^{1/r}.$$

We get that the norm of the bilinear operator $S_{\alpha,\rho}$ satisfies
\begin{eqnarray*}
\|S_{\alpha,\rho}\|_{p,q,r} \geq N^{1/r-1/q-1/2}= N^{1/p-1/2}.
\end{eqnarray*}
A similar argument with the roles of $f$ and $g$ interchanged yields 
\begin{eqnarray*}
\|S_{\alpha,\rho}\|_{p,q,r} \geq N^{1/q-1/2}.
\end{eqnarray*}
This completes the proof. 
 
%\end{proof}
\qed
\subsection{Proof of Lemma~\ref{ml}}\label{mlp}

For a cube $Q$ in $\R^n$ and $\alpha>0$, let $\alpha Q$ denote the cube with same center as of $Q$ and measure $\alpha^n |Q|$, where $|Q|$ denote the $n$ dimensional Lebesgue measure of $Q$. We denote by $l(Q)$ the sidelength of the cube $Q$. 

We know that bounded bilinear operators have $l_2-$valued extension, i.e., if $T$ is a bounded bilinear operator from $L^p(\R^n)\times L^q(\R^n)$ into $L^r(\R^n)$, then for finite sequences of functions $\{f_k\}_{k=1}^N$ and $\{g_k\}_{k=1}^N$, we have 
\begin{eqnarray}\label{vve}
\|\{T(f_k,g_k)\}\|_{L^r(l_2)} \lesssim \|\{f_k\}\|_{L^p(l_2)}\|\{g_k\}\|_{L^q(l_2)},
\end{eqnarray}
where the implicit constant is a multiple of the operator norm $\|T\|_{p,q,r}$. See Section 7 in~\cite{gm} for vector-valued extension of bilinear operators. 
Let $Q_k$'s be as defined in the statement of Lemma $2.2$.
Let $\{u_k\}_{k=1}^N$ and $\{v_k\}_{k=1}^N$ be finite sequences of compactly supported smooth functions. Fix an unbounded sequence of real numbers $\{\lambda_m \}_{m\in\Z}$. Take $t_k\in \frac{1}{2} Q_k,~1\leq k\leq N$ and set $f_k^m(t)=u_k(\lambda_m (t-t_k))$ and $g_k^m(t)=v_k(\lambda_m t)$. 
Since $\lambda_m\rightarrow \infty$, we can choose $\lambda_m$ large enough to ensure that 
$A_k=\text{supp}(f_k^m)\subset \frac{1}{2} Q_k$ and $B_k=\text{supp}(g_k^m)\subset [\frac{-l(Q_k)}{4},  \frac{l(Q_k)}{4}]^n.$ Observe that 
$A_k-B_k:=\{x-y: x\in A_k, y\in B_k\} \subset Q_k$ for all $k$. 
%Define $F_k=\check f_k$ and $G_k=\check g_k$. 

Let $T_{\phi,\lambda_m}$ be the bilinear operator associated with the function $e^{i\lambda_m \phi(\xi-\eta)}$, i.e., 
\begin{eqnarray*}
T_{\phi,\lambda_m}(f, g)(x)
%&=& \int_{\R^n} \int_{\R^n} \hat f_k(\xi)\hat G_k(\eta) e^{i\lambda_m \phi(\xi-\eta)} e^{2\pi i x\cdot(\xi+\eta)} d\xi d\eta\\
&=& \int_{\R^n} \int_{\R^n} \hat f(\xi)\hat g(\eta) e^{i\lambda_m \phi  (\xi-\eta)} e^{2\pi i x\cdot(\xi+\eta)} d\xi d\eta.
\end{eqnarray*}
Since $A_k-B_k\subset Q_k$ and $\phi(t)=\langle \alpha_k, t\rangle +\beta_k$ for almost all $t\in Q_k$, we get that 
\begin{eqnarray*}
|T_{\phi,\lambda_m}(\check f_k^m,\check g_k^m)(x)|=\left | \check f_k^m(x+\lambda_m \alpha_k) \check g_k^m(x-\lambda_m \alpha_k)\right |,
\end{eqnarray*}
where $\check f(x)=\int_{\R^n} f(y) e^{2\pi ix.y} dy.$

Further, note that $ |\check f_k^m(t)|=|\frac{1}{\lambda_m^n}\check u_k(\frac{t}{\lambda_m})|$ and $ |\check g_k^m(t)|=|\frac{1}{\lambda_m^n}\check v_k(\frac{t}{\lambda_m})|.$

Therefore, using the vector-valued extension~(\ref{vve}) for the operator $T_{\phi,\lambda_m}$ along  with a change of variables, we get
\begin{eqnarray*}
\left \|\left( \sum_{k=1}^N| \check u_k(x+\alpha_k) \check v_k(x-\alpha_k)|^2\right)^{\frac{1}{2}}\right \|_{L^r_x} &=&
 \lambda^{2n}_m \left \|\left( \sum_{k=1}^N| \check f_k^m(\lambda_m(x+\alpha_k))\check g_k^m(\lambda_m(x-\alpha_k))|^2\right)^{\frac{1}{2}}\right \|_{L^r_x} \\
 &=&
 \lambda^{(2-1/r)n}_m \left \|\left( \sum_{k=1}^N| \check f_k^m(x+\lambda_m\alpha_k)\check g_k^m(x-\lambda_m\alpha_k)|^2\right)^{\frac{1}{2}}\right \|_{L^r_x} \\
% &=& \left \|\left( \sum_{k=1}^N \left |\frac{1}{\lambda_m} F_k(\frac{x}{\lambda_m}+\alpha_k) \frac{1}{\lambda_m}(\frac{x}{\lambda_m}-\alpha_k)\right | ^2\right)^{\frac{1}{2}}\right \|_{L^r_x} \\
 &=&  \lambda^{(2-1/r)n}_m\left \| \left( \sum_{k=1}^N | T_{\phi,\lambda_m}(\check f_k^m,\check g_k^m)(x)| ^2\right)^{\frac{1}{2}}\right \|_{L^r_x}\\
 &\lesssim &   \lambda^{-n/r}_m\left \| \left( \sum_{k=1}^N |\check u_k(\frac{x}{\lambda_m})| ^2\right)^{\frac{1}{2}}\right \|_{L^p_x}  \left \| \left( \sum_{k=1}^N | \check v_k(\frac{x}{\lambda_m})| ^2\right)^{\frac{1}{2}}\right \|_{L^q_x}\\
 &= & C_{\phi} \left \| \left( \sum_{k=1}^N |\check u_k| ^2\right)^{\frac{1}{2}}\right \|_{L^p}  \left \| \left( \sum_{k=1}^N | \check v_k| ^2\right)^{\frac{1}{2}}\right \|_{L^q},
\end{eqnarray*}
where the constant $C_{\phi}$ is a constant multiple of $\sup_{k\in \mathbb Z} \|e^{\i\lambda_k \phi(\xi-\eta)}\|_{\mathcal M_{p,q,r}(\R^n)}$.

Since functions with compactly supported Fourier transform are dense in $L^p$ spaces, $1\leq p<\infty$, we see that the vector-valued operator $S_{\alpha,1}$ for $\alpha=(\alpha_1,\alpha_2,\dots,\alpha_N)$, extends boundedly from $L^p(l_2)\times L^q(l_2)$ to $L^r(l_2)$. 

This together with Proposition~\ref{vector} yields 
$$\sup_{k\in \mathbb Z} \|e^{i\lambda_k \phi(\xi-\eta)}\|_{\mathcal M_{p,q,r}(\R^n)}\geq \max\{N^{1/p-1/2}, N^{1/q-1/2}\}.$$

Next, by considering an adjoint (as defined in Section~\ref{bil}) of the bilinear multiplier operator $T_{\phi,\lambda_m}$, a similar argument as above can be used to prove the final estimate ~(\ref{norm}) for the operator norm of bilinear multipliers under consideration. 

This completes the proof of Lemma~\ref{ml}.
%Further, the number $N$ is arbitrary, therefore Proposition~\ref{vector} would lead to a contradiction if $p<2$ or $q<2$. Further, using the transpose of bilinear operator $T_m$ we can conclude the same when $r>2$. \color{red} {THIS NEEDS TO BE VERIFIED.}

%%%%%%%%%%%%%%%%%%%%%%%%%%%%%%%%

\subsection{Proof of Theorem~\ref{mt}}
The proof of the theorem is by contradiction. Let us assume on the contrary that for a sequence of positive numbers $\{\lambda_m\}_{\N}$ with $\lambda_m \rightarrow \infty$ as $m\rightarrow \infty,$ we have the uniform control over multiplier norm, i.e., 
$$\|e^{i\lambda_m \phi(\xi-\eta)}\|_{\mathcal M_{p,q,r}(\R^n)}\leq C,~ ~~\forall m.$$

Let $Q\subset \R^n$ be a cube such that $\phi$ is non-linear on $Q$. Then for a fixed $N\in \N$, we can choose distinct points $t_k,~1\leq k\leq N,$ in $Q$ such that the gradient $\nabla \phi(t_k)=\alpha_k$ are distinct for all $1\leq k\leq N$. 

Define $s_k(t)=\phi(t_k)+\alpha_k\cdot (t-t_k).$ 

Let $\{u_k\}_{k=1}^N$ and $\{v_k\}_{k=1}^N$ be finite sequences of compactly supported smooth functions. Choose $\Lambda>0$ so that 
$\text{supp}(u_k), \text{supp}(v_k)\subset [-\Lambda,\Lambda]^n.$ 

Consider $f^m_k(t)=u_k(\lambda_m (t-t_k))$ and $g^m_k(t)=v_k(\lambda_m t)$. 

Since $\lambda_m\rightarrow \infty$, we can choose $\lambda_m$ large enough to ensure that 
$A_k=\text{supp}(f_k^m)\subset \frac{1}{2} Q_k$ and $B_k=\text{supp}(g_k^m)\subset [\frac{-l(Q_k)}{2},  \frac{l(Q_k)}{2}]^n,$ where $Q_k=[t_k-\frac{\Lambda}{\lambda_m},  t_k+\frac{\Lambda}{\lambda_m}]^n.$ Observe that 
$A_k-B_k:=\{x-y: x\in A_k, y\in B_k\} \subset Q_k$ for all $k$. 
%Define $F_k=\check f_k$ and $G_k=\check g_k$. 
Since $\phi\in C^1(\R^n)$, we have 
$$\sup\limits_{t\in Q_k}|\phi(t)-s_k(t)|=o\left(\frac{1}{\lambda_m^n}\right),~m\rightarrow \infty.$$
Let $T$ and $T_k$ be the bilinear multiplier operators associated with functions $e^{i\lambda_m \phi(\xi-\eta)}$ and $e^{i\lambda_m s_k(\xi-\eta)}$ respectively. Note that 
\begin{eqnarray*}
|T(\check f_k^m, \check g_k^m)(x)-T_k(\check f_k^m, \check g_k^m)(x)|
&\leq & \int_{\R^n}\int_{\R^n} |f_k^m(\xi)||g_k^m(\eta)| |e^{i\lambda_m\phi(\xi-\eta)}-e^{i\lambda_m s_k(\xi-\eta)}|d\xi d\eta\\
&=& o\left(\frac{1}{\lambda_m^{3n-1}}\right)~m\rightarrow \infty.
\end{eqnarray*}
Therefore, uniformly in $x\in \R^n$, we get the following 
\begin{eqnarray}
\left(\sum\limits_{k=1}^N |T_k(\check f_k^m, \check g_k^m)(x)|^2\right)^{\frac{1}{2}}
&\leq & \left(\sum\limits_{k=1}^N |T(\check f_k^m, \check g_k^m)(x)|^2\right)^{\frac{1}{2}}+o\left(\frac{1}{\lambda_m^{3n-1}}\right),
\end{eqnarray}
as $m\rightarrow \infty.$
Now for a fixed $\Gamma>0$, 
\begin{eqnarray*}
\left(\int_{|x|<\Gamma \lambda_m} \left(\sum\limits_{k=1}^N |T_k(\check f_k^m, \check g_k^m)(x)|^2\right)^{\frac{r}{2}}\right)^{\frac{1}{r}}
&\leq & \left(\int_{|x|<\Gamma \lambda_m} \left(\sum\limits_{k=1}^N |T(\check f_k^m, \check g_k^m)(x)|^2\right)^{\frac{r}{2}}\right)^{\frac{1}{r}}\\ 
&&+o\left(\frac{1}{\lambda_m^{3n-1-\frac{n}{r}}}\right),~~m\rightarrow \infty\\
&\leq & \left(\int_{\R^n} \left(\sum\limits_{k=1}^N |\check f_k^m|^2\right)^{\frac{p}{2}}dx\right)^{\frac{1}{p}} \left(\int_{\R^n} \left(\sum\limits_{k=1}^N |\check g_k^m|^2\right)^{\frac{q}{2}}dx\right)^{\frac{1}{q}}\\ 
& &+o\left(\frac{1}{\lambda_m^{3n-1-\frac{n}{r}}}\right),~~m\rightarrow \infty
\end{eqnarray*}
Using the definitions of $T_k, f_k^m,$ and $g_k^m$, we get that 
$
|T_k(\check f_k^m,\check g_k^m)(x)|=\left | \check f_k^m(x+\lambda_m \alpha_k) \check g_k^m(x-\lambda_m \alpha_k)\right |$,  $|\check f_k^m(t)|=|\frac{1}{\lambda_m^n}\check u_k(\frac{t}{\lambda_m})|$ and $ |\check g_k^m(t)|=|\frac{1}{\lambda_m^n}\check v_k(\frac{t}{\lambda_m})|.$
After substituting these in the inequality above and simplifying it further, we finally get the following 

\begin{eqnarray*}
\left (\int_{|x|<\Gamma} \left( \sum_{k=1}^N| \check u_k(x+\alpha_k) \check v_k(x-\alpha_k)|^2\right)^{\frac{r}{2}}dx\right)^{\frac{1}{r}}
&\lesssim &  \left \| \left( \sum_{k=1}^N |\check u_k| ^2\right)^{\frac{1}{2}}\right \|_{L^p}  \left \| \left( \sum_{k=1}^N | \check v_k| ^2\right)^{\frac{1}{2}}\right \|_{L^q}\\
&& +o\left(\frac{1}{\lambda_m^{n-1}}\right),~~m\rightarrow \infty. 
\end{eqnarray*}
Since $\Gamma>0$ is arbitrary, we get a contradiction using Lemma~\ref{ml}. 
\begin{remark}
If $\phi$ is a homogeneous function of degree at least $2$, then using the dilation invariance property of bilinear multipliers, it can be deduced from our theorem that it cannot be a bilinear multiplier for exponents lying outside the local $L^2-$range. An important example of such functions is $\phi(\xi)=|\xi|^2$. 
\end{remark}
\section{Counterexamples in Local $L^2$ region}
\subsection{The function $e^{i|\xi-\eta|^2}$ in local $L^2-$ range} In this section we discuss examples of non-linear functions $\phi$ for which $e^{i \phi(\xi-\eta)}$ does not give rise to bilinear multiplier for exponents in the local $L^2-$range. Therefore, even in local $L^2-$ range, we cannot expect to have a consistent positive result concerning bilinear multipliers of the form $e^{i\phi(\xi-\eta})$, where $\phi$ is a non-linear function. We also discuss some other examples of unimodular bilinear multipliers for which the boundedness fails in the local $L^2-$ range.

\begin{center}
\includegraphics[scale=0.5]{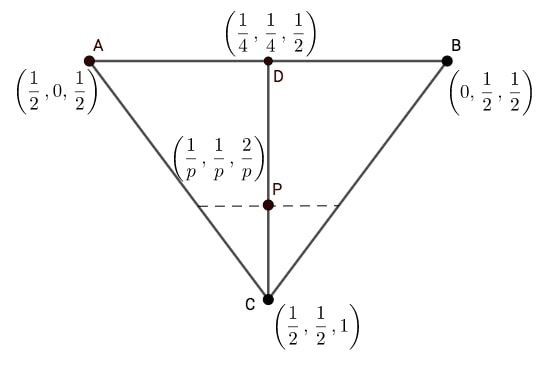}\\
{Figure 1:~Local $L^2-$ range}
\end{center}
Note that a triplet $(p,q,r)$ is in the local $L^2-$ range whenever the point $\left(\frac{1}{p},\frac{1}{q},\frac{1}{r}\right)$ is in the region (including the boundary) as described in Figure 1.

We first consider $\phi(\xi-\eta)= |\xi-\eta|^2$ and show that this does not give rise to unimodular bilinear multiplier for any triplet in the local $L^2-$ range.
Let $$T(f,g)(x):= \int_{\R^n}\int_{\R^n}e^{-i|\xi-\eta|^2}\hat f(\xi)\hat g(\eta) e^{i x.(\xi+\eta)}d\xi d\eta.$$ 
Recall that if $T$
is bounded for a fixed triplet $(p,q,r)$ satisfying the H\"{o}lder condition, then by the standard dilation argument the bilinear operator $T_t,$ associated with symbol $e^{it |\xi-\eta|^2}$, is also bounded for the same triplet with $\|T\|_{\mathcal{M}_{p,q,r}(\R^n)}=\|T_t\|_{\mathcal{M}_{p,q,r}(\R^n)}.$\\

For fixed $f,g,h$ in the Schwartz class $\mathcal{S}(\R^n)$, we write 
\begin{eqnarray*}
\langle f g, h\rangle &=&  \int_{\R^n}\int_{\R^n} \hat f(\xi)\hat g(\eta) \hat h (\xi+\eta) d\xi d\eta \\
&=& \int_{\R^n}\int_{\R^n} \hat f(\xi)\hat g(\eta)  e^{-it|\xi-\eta|^2} e^{it|\xi-\eta|^2} e^{it|\xi+\eta|^2}e^{-it|\xi+\eta|^2}\hat h (\xi+\eta) d\xi d\eta \\
&=& \int_{\R^n}\int_{\R^n}  e^{i2t|\xi|^2}\hat f(\xi) e^{i2t|\eta|^2}\hat g(\eta) e^{-it|\xi-\eta|^2} e^{-it|\xi+\eta|^2} \hat h (\xi+\eta) d\xi d\eta \\
&=&\langle T_t(Q_{2t}f,Q_{2t}g),Q_{-t}h\rangle,
\end{eqnarray*}
where $Q_t$ denotes the Fourier multiplier operator $\widehat{Q_t f}(\xi)=e^{it|\xi|^2}\hat f(\xi).$ 
In \cite{h}, L. Hormander proved that $\|Q_{t}f\|_p \leq |t|^{\frac{1}{p}-\frac{1}{2}}\|f\|_2$ for $p\geq 2$ and $t\neq 0.$ If the bilinear operator $T_t$ is bounded for a triplet $(p,q,r)$ in local $L^2-$ range then using the expression above with the estimate on $Q_t f$, we get that 
\begin{eqnarray*}
|\langle f g, h\rangle | &\leq & \|T\|_{\mathcal{M}_{p,q,r}(\R^n)} \|Q_{2t} f\|_p \|Q_{2t} g\|_q \|Q_{-t}h\|_{r'}\\
&\leq & C(p,q,r) \|T\|_{\mathcal{M}_{p,q,r}(\R^n)} |t|^{\frac{1}{p}+\frac{1}{q}+\frac{1}{r'}-\frac{3}{2}}\\
&=& C(p,q,r) \|T\|_{\mathcal{M}_{p,q,r}(\R^n)} |t|^{-\frac{1}{2}}.
\end{eqnarray*} 
Letting $|t| \rightarrow \infty$ we get a contradiction. Thus, $T$ cannot be a bounded operator for triplets lying in local $L^2-$ range. 
\subsection{The function $e^{i\xi\cdot \eta}$ on the boundary of local $L^2-$ range}
We observe that the function $e^{i\xi\cdot \eta}$ cannot be a bilinear multiplier for points on the boundary which consists of line segments $AC, CB$ and $BA$ (see Figure 1). 

First, verify that the boundedness of the bilinear operator $T$ associated with $e^{i\xi\cdot \eta}$ at one end-point (vertex of triangle) implies the same at the other two end-points. For example, if $T$ bounded from $L^2\times L^2\rightarrow L^1$, then by considering adjoints  $T^{*,1}$ and $T^{*,2}$, the multiplier takes the form $e^{-i(\xi+\eta).\eta}$ and $e^{-i\xi.(\xi+\eta)}.$  So it can be easily deduced that $T$ is also bounded from $L^{\infty} \times L^2 \rightarrow L^2$ and $L^2\times L^{\infty}\rightarrow L^2$ respectively. 

Since the function $e^{i\xi\cdot \eta}$ is symmetric in $\xi$ and $\eta$, it is enough to show that the operator $T$ is not bounded from $L^2\times L^{\infty}\rightarrow L^2$. Consider,
\begin{eqnarray*}\label{123}
\langle T(f,g), h\rangle &=&  \int_{\R^n}\int_{\R^n} \hat f(\xi)\hat g(\eta) e^{2i\xi\cdot \eta} \hat h (\xi+\eta) d\xi d\eta \\
&=& \int_{\R^n}\int_{\R^n} \hat f(\xi)\hat g(\eta) e^{i(|\xi+\eta|^2-|\xi|^2|-\eta|^2)} \hat h (\xi+\eta) d\xi d\eta \\
&=& \langle L_{\bar {\sigma}}(f) L_{\bar {\sigma}}(g), L_{\sigma} h\rangle,
\end{eqnarray*}
where $\widehat {L_{\bar {\sigma}}(f)}= e^{-i|\xi|^2} \hat f$ and $\widehat {L_{\sigma}(f)}= e^{i|\xi|^2} \hat f$.

Note that the relation above holds for all $f,h \in L^2$ and $g\in L^{\infty}$. Therefore, we get that 
\begin{eqnarray*}
\left| \langle f L_{\bar {\sigma}}(g), h\rangle \right| &=& \left|\langle T(L_{\sigma}(f),g), L_{\bar {\sigma}}(h)\rangle \right|\\
&\lesssim & \|f\|_2 \|g\|_{\infty} \|h\|_2.
\end{eqnarray*}
The duality implies that 
$$\|f L_{\bar {\sigma}}(g)\|_2\ \lesssim \|f\|_2 \|g\|_{\infty}$$
for all $f\in L^2$ and $g\in L^{\infty}$. From here we conclude that 
$L_{\bar {\sigma}}$ is bounded from $L^{\infty}$ into itself and this leads to a contradiction. 

Next, let us consider points on the line segment $BA$. Note that for a given point $(\frac{1}{p},\frac{1}{q},\frac{1}{2})$ on the line segment $BA$, if the function $e^{2i\xi\cdot \eta}\in \mathcal{M}_{p,q,2}(\R^n)$, then due to the symmetry of $e^{2i\xi\cdot \eta}$ with respect to $\xi$ and $\eta$, we always have that $e^{2i\xi\cdot \eta}\in \mathcal{M}_{q,p,2}(\R^n)$. Therefore, by the interpolation argument, it  is enough to show that $e^{2i\xi\cdot \eta}$ is not  a bilinear multiplier at the point $D=(\frac{1}{4},\frac{1}{4},\frac{1}{2})$ (see Figure 1). 

From the calculations as above, we have
\begin{eqnarray*}
\langle T(f,f), L_{\bar {\sigma}} h\rangle 
&=& \langle L_{\bar {\sigma}}(f)^2 , h\rangle.
\end{eqnarray*} 
Thus, the duality and the boundedness of $T$ from $L^4\times L^4$ to $L^2$ would imply the boundedness of the linear operator $L_{\bar {\sigma}}$ on $L^4,$
which is not possible.  

For the remaining two line segments $AC$ and $CB$, again note that due to the symmetry, it is enough to consider the case of one line segment. Consider a point $(\frac{1}{2},\frac{1}{q},\frac{1}{r})$ on $AC$. If $e^{2i\xi\cdot \eta}\in \mathcal{M}_{2,q,r}(\R^n)$, then by considering the adjoint we get that $e^{-2i(\xi\cdot \eta+|\eta|^2)}\in \mathcal{M}_{r',q,2}(\R^n)$.
Consider the bilinear multiplier operator associated with the multiplier symbol $e^{-2i(\xi\cdot \eta+|\eta|^2)}$ and write
\begin{eqnarray*}
\langle T(f,g), h\rangle &=&  \int_{\R^n}\int_{\R^n} \hat f(\xi)\hat g(\eta) e^{-2i(\xi\cdot \eta+|\eta|^2)} \hat h (\xi+\eta) d\xi d\eta \\
&=& \int_{\R^n}\int_{\R^n} \hat f(\xi)\hat g(\eta)e^{-2i(\xi\cdot \eta+|\eta|^2)}  e^{i|\xi+\eta|^2} e^{-i|\xi+\eta|^2} \hat h (\xi+\eta) d\xi d\eta \\
&=& 
\langle L_{\sigma}(f) L_{\bar {\sigma}}(g), L_{\bar{\sigma}} h\rangle.
\end{eqnarray*}
As $L_{\bar{\sigma}} $ is a unitary operator on $L^2$, if $T$ is a bilinear multiplier for triplet $(r',q,2)$ then the symbol corresponding to $L_{\sigma}(\cdot) L_{\bar {\sigma}}(\cdot)$ belongs to $\mathcal{M}_{r',q,2}(\R^n).$
An easy computation shows that $\overline{L_{\sigma}f} = L_{\bar{\sigma}}\bar{f}.$ This immediately gives us that symbol of $L_{\sigma}(\cdot) L_{\bar {\sigma}}(\cdot)$ belongs to $\mathcal{M}_{q,r',2}(\R^n).$ So it is enough to look at the case $r'=q=4.$
 At this point a similar argument, as in the previous case, can be used to get a contradiction.

\subsection{The function $e^{i(|\xi|^2+|\eta|^2)}$ in the interior of local $L^2-$ range} This case follows in a similar fashion. We claim that the function $e^{i(|\xi|^2+|\eta|^2)}$ does not give rise to a bilinear multiplier for any point in the interior of the region described in Figure 1. As previously, using the symmetry of the function $e^{i(|\xi|^2+|\eta|^2)},$ it is enough to show that $e^{i(|\xi|^2+|\eta|^2)}\notin \mathcal{M}_{p,p,\frac{p}{2}}(\R^n)$ for $2<p<4$.  For, if $e^{i(|\xi|^2+|\eta|^2)}\in \mathcal{M}_{p,p,\frac{p}{2}}(\R^n)$ with $2<p<4$ then using exactly the same argument as in the previous case one can show that the linear operator $L_{{\sigma}}$ is bounded on $L^p,$ which is not the case.

\begin{remark}We would like to point out that the boundedness properties of bilinear multiplier operators associated with symbols $e^{i|\xi-\eta|^2}$ and $e^{i\xi\cdot \eta}$ are equivalent at $(2,2,1), (2,\infty, 2)$ and $(\infty, 2,2)$. This can be verified by the relation $e^{i|\xi-\eta|^2}=e^{i|\xi|^2+i|\eta|^2-2i\xi\cdot \eta}$ and using the adjoints of bilinear operator. Therefore, the function $e^{i|\xi-\eta|^2}$ cannot give rise to bilinear multiplier at any of the three end-points of the local $L^2-$triangle.   
\end{remark}

\section{Essential continuity}\label{ec} 
In this section we study the essential continuity property of unimodular bilinear multipliers. Let us first define the term essential continuity of a function. 
\begin{definition}Let $f$ be a measurable function in $\R^n$. A point $x\in \R^n$ is said to be a point of essential continuity of $f$ if for every $\epsilon>0$ there is a neighbourhood $B_x$ of $x$ such that 
$|f(x)-f(y)|< \epsilon~~\text{for~ almost~ every~} y\in B_x.$

Let  $\Omega(\psi,x):= \lim_{\delta\rightarrow 0} \sup_{|x-y|<\delta}|f(x)-f(y)|.$

\end{definition}
\begin{definition}
Let $E$ be a measurable set in $\R^n.$ Then a point $x\in\R^n$ is called a density point of $E$ if $$\lim_{\delta\rightarrow 0}\frac{|E\cap B(x,\delta)|}{|B(x,\delta)|} =0.$$
We denote the set of all density points of $E$ by $E^d.$ It is known that the symmetric difference  $E\triangle E_d$ has Lebesgue measure zero.

\end{definition}
\begin{theorem} Let $(p,q,r)$ be a triplet lying outside local $L^2-$ range and $\psi$ be a measurable function such that $\psi(\xi-\eta)$ and $|\psi(\xi-\eta)|^2$ defines a bounded bilinear multiplier from $L^p(\R^n)\times L^q(\R^n)$ to $L^r(\R^n),$  then the function $\psi(x)$ is essentially continuous at almost every point $x\in \R^n$.
\end{theorem}  
%\begin{proof}
{\bf Proof:} 
The proof is by contradiction. 
 Let $\Omega(\psi,x)$ be the essential oscillation of $\psi$ at $x\in \R^n$. 
Suppose on the contrary that the set of all points at which $\psi$ is not essential continuous is of positive measure. 

The above assumption on $\psi$ implies that there is an $\epsilon>0$ so that the set 
$$E:=\{x\in \mathcal L(\psi): \Omega(\psi,x)>\epsilon\}$$ 
has positive measure. 

For a real number $c$ consider the sets
$$E_1:=\{x \in E: |\psi(x)-c|<\frac{\epsilon}{3}\}~\text{and}~E_2:=\{x\in \mathcal L(\psi): |\psi(x)-c|>\frac{2\epsilon}{3}\}$$ 

It is easy to see that for some scalar $c$, the set $E_1$ has positive measure using the density of rational numbers.  Moreover 
$$|\bar{E^{d}_1}\cap \bar{E_2^d}|>0.$$
Consider $m(\xi):=9|\psi(\xi)-c|^2-3\epsilon^2$. By our assumption, $|\psi(\xi-\eta)|^2$ and $\psi(\xi-\eta)$ define a bounded bilinear multiplier from $L^p(\R^n)\times L^q(\R^n)$ to $L^r(\R^n),$  so on expanding the expression of $m$ above we see that $m(\xi-\eta)\in \mathcal M_{p,q,r}(\R^n)$ and for $\delta>0$, the function $m_{\delta}:=m\ast \frac{1}{|B(0,\delta)|} \chi_{_{B(0,\delta)}}$ is also a bilinear multiplier for the same tuple $(p,q,r)$. 

Note that $m(\xi)<-\epsilon^2$ for $\xi\in E_1$ and 
$m(\xi)>\epsilon^2$ for $\xi\in E_2.$

For a positive integer $N$ consider a sequence $\{c_k\}_{k=1}^N$ with $c_k\in \{-1,1\}$. Lemma 1 (Lebedev and Olevskii \cite{o2} p. 554) guarantees an arithmetic progression $t_k=a+kh,~k=1,2,3,\dots,N$ with the properties that $t_k\in E^d_1$ if $c_k=-1$ and $t_k\in E^d_2$ if $c_k=1$. 

Use the transference principle for bilinear multipliers and deduce that the restriction of $m_{\delta}$ to the arithmetic progression give rise to a bilinear multiplier operator, i.e, 
$m_{\delta}(t_k)\in \mathcal M_{p,q,r}(\Z^n)$. Moreover, the norm of the multiplier operator is bounded uniformly in $\delta$. The restriction to $t_k$ makes sense as these points are Lebesgue points of $m$. Since $m_{\delta}(t_k)\rightarrow m(t_k)$ as $\delta \rightarrow 0$, we know that $m(t_k)\in \mathcal M_{p,q,r}(\Z^n)$ as well.

For the trigonometric polynomial $P(x)=\sum\limits_{k=1}^Nc_ke^{2\pi i k\cdot x}$ consider 
\begin{eqnarray*}
\sum\limits_{k=1}^N c_km_{\delta}(t_k)&=& \sum\limits_{k=1}^N \left(\int_{\T^n} P(x)e^{2\pi i k\cdot x} dx\right)m_{\delta}(t_k)\\
&=& \int_{\T^n}\left(\sum\limits_{k=1}^N m_{\delta}(t_k) e^{2\pi i k\cdot x} \right)P(x)dx\\
&=& \int_{\T^n} M_{m_{\delta}}(Q,{\bf 1})P(x)dx
\end{eqnarray*}
where $Q(x)=\sum\limits_{k=1}^N e^{2\pi i k\cdot x}, {\bf 1}=\chi_{[0,1)^n}$ and $M_{m_{\delta}}$ is the bilinear multiplier operator associated with symbol $m_{\delta}$ on the torus group.
Using the boundedness of the operator $M_{m_{\delta}}$ together with the choice of $t_k$, we get 
\begin{eqnarray*}
\epsilon^2 N=\sum\limits_{k=1}^N c_km_{\delta}(t_k)&=&|\int_{\T^n}M_{m_{\delta}}(Q,{\bf 1})P(x)dx|\\
&\leq & \|M_{m_{\delta}}(Q,{\bf 1})\|_r \|P\|_{r'}\\
&\leq & \|Q\|_p \|P\|_{r'}.
\end{eqnarray*}
Since the tuple $(p,q,r)$ is outside the local $L^2-$range. There are three possibilities: $p<2$ or $q<2$ or $r'<2$. 

Let us first consider the case when $p<2$. This means $r<2~\text{or}~r'>2$. In this case we can choose $P$, in the beginning itself, so that $\|P\|_{r'}\leq N^{1/2}.$ Further, using inclusion relation of $L^p(\T^n)$ spaces we have $\|Q\|_p\lesssim N^{1/p'}.$
Putting everything together we get 
$$N\lesssim N^{1/2+1/p'}.$$
This is possible only when 
$1/2+1/p'\geq 1~\text{or}~p\geq 2,$ which is a contradiction. 

The other case when $q<2$ can be dealt with similarly using the symmetry. 

Finally, when $r'<2$. In this case $p,q,r>2$. We shall consider the first transpose of operator $M_{m_{\delta}}$. This will reduce the problem to the previously known case.

This completes the proof.
\qed
%\end{proof}

We get the following corollary:
\begin{corollary} Let $\phi:\R^n\rightarrow [0,2\pi)$ be a measurable function. If $e^{i\lambda \phi(\xi-\eta)} \in \mathcal M_{p,q,r}(\R^n)$, where $ (p,q,r)$ lies outside the local $L^2-$ range with a uniform bound on the multiplier norms
 $$\|e^{i\lambda \phi(\xi-\eta)}\|_{\mathcal M_{p,q,r}(\R^n)}\leq C,\forall \lambda\in \Z,$$ then the function $\phi$ is essential continuous at almost every point $x\in \R^n$.  
\end{corollary}
%\begin{proof}
{\bf Proof:}
Write $\psi(\xi-\eta)=e^{i\phi(\xi-\eta)}$. Note that $|\psi(\xi-\eta)|^2 =1$, it trivially defines a bounded bilinear multiplier for any triplet $(p.q.r).$ Therefore, for exponents $ (p,q,r)$ outside the local $L^2-$ range, the last theorem implies that $e^{i\phi}$ is an essentially continuous function. 

Further, the essential continuity of $\phi$ may be proved imitating the arguments given in Lemma 2 and Lemma 3 in~\cite{o2}. Here we would also  require the result about idempotent multipliers in bilinear setting from ~\cite{S}. Since, the proof of this assertion follows verbatim, we skip the details here. 
\qed

\section*{Acknowledgement} The authors would like to thank the referee for valuable suggestions which greatly helped us in improving the presentation of the results. The second author acknowledges the financial support from the Department of Science and Technology, Government of India under the scheme MATRICS with project MAT/2017/000039/MS.

\end{document}